\documentclass[12pt]{amsart}
\usepackage{amssymb}
\usepackage{amsfonts}
\usepackage{latexsym}
\usepackage{amscd}
\usepackage[mathscr]{euscript}
\usepackage{xy} \xyoption{all}

\newcommand{\uloopr}[1]{\ar@'{@+{[0,0]+(-4,5)}@+{[0,0]+(0,10)}@+{[0,0] +(4,5)}}^{#1}}
\newcommand{\uloopd}[1]{\ar@'{@+{[0,0]+(5,4)}@+{[0,0]+(10,0)}@+{[0,0]+ (5,-4)}}^{#1}}
\newcommand{\dloopr}[1]{\ar@'{@+{[0,0]+(-4,-5)}@+{[0,0]+(0,-10)}@+{[0, 0]+(4,-5)}}_{#1}}
\newcommand{\dloopd}[1]{\ar@'{@+{[0,0]+(-5,4)}@+{[0,0]+(-10,0)}@+{[0,0 ]+(-5,-4)}}_{#1}}

\newcommand{\luloop}[1]{\ar@'{@+{[0,0]+(-8,2)}@+{[0,0]+(-10,10)}@+{[0, 0]+(2,2)}}^{#1}}

\newtheorem{lem}{Lemma}[section]
\newtheorem{corol}[lem]{Corollary}
\newtheorem{theor}[lem]{Theorem}
\newtheorem{prop}[lem]{Proposition}
\theoremstyle{definition}
\newtheorem{defi}[lem]{Definition}

\newtheorem{exem}[lem]{Example}

\newtheorem{rema}[lem]{Remark}

\newcommand{\Qs}{Q_{\sigma}}
\newcommand{\Dl}{\mathcal{I}_{dl}}
\newcommand{\Dr}{\mathcal{I}_{dr}}

\begin{document}

\title[Two-sided localization of bimodules]{Two-sided localization of bimodules}%
\author{E. Ortega}
\address{Departament de Matem\`atiques, Universitat Aut\`onoma de Barcelona,
08193 Bellaterra (Barcelona), Spain.} \email{eortega@mat.uab.es}

\thanks{} 
\date{\today}

\begin{abstract} We extend to bimodules Schelter's localization of a
ring with respect to Gabriel filters of left and right ideals. Our
two-sided localization of bimodules provides an endofunctor on a
convenient bicategory of rings with filters of ideals. This is used
to study the Picard group of a ring relative to a filter of ideals.
\end{abstract}
\maketitle
\section*{Introduction}
Given a right $R$-module $M_R$ and a certain filter $\Im$ of right
ideals of $R$, Gabriel \cite{Gabriel} and Goldman \cite{Goldman}
gave the construction of the module of quotients of $M$ with respect
to the filter $\Im$, denoted by $Q_{\Im}(M)$. Moreover, if we
consider the regular right $R$-module $R$ then $Q_{\Im}(R)$ turns
out to be a ring. It is shown  in \cite{Bo} that the more classical
rings of quotients, as for example the Ore localizations, are
particular cases of this construction. Moreover, this localization
produces a functor $M\mapsto Q_{\Im}(M)$ from the category of right
$R$-modules to the category of right $Q_{\Im}(R)$-modules. In
\cite{Schelter2} it is considered the category $\mathcal{W}$ of
rings with a Gabriel filter of right ideals, and the endofunctor
$Q(R,\Im)=(Q_{\Im}(R),\Im')$. In case that $R$ is a commutative
ring, it is shown that the subcategory $\mathcal{C}$ of the complete
pairs, i.e. those such that $(R,\Im)\cong (Q_{\Im}(R),\Im')$, is
reflective in $\mathcal{W}$. This is not true for general
noncommutative rings.

Later on, Schelter \cite{Schelter} generalized the commutative case
by considering certain Gabriel filters on $R\otimes_{\mathbb{Z}}
R^{op}$. In particular, he considered the category
$\mathcal{W}_{\sigma}$ of rings $R$ with Gabriel filters
$\mathcal{D}_r$ and $\mathcal{D}_l$ of right and left ideals
respectively. He constructed the two-sided localization of $R$ with
respect to $\mathcal{D}_r$ and $\mathcal{D}_l$, and obtained a ring
structure $Q_{\mathcal{D}_r,\mathcal{D}_l}(R)$. In the particular
case where $\Dr$ and $\Dl$ are the Gabriel filters of right and left
dense ideals respectively, $Q_{\Dr,\Dl}(R)=\Qs(R)$ is the maximal
symmetric ring of quotients (see \cite{Lanning} and \cite{Ortega}).
This construction defines an endofunctor
$(R,\mathcal{D}_l,\mathcal{D}_r)\mapsto
(Q_{\mathcal{D}_r,\mathcal{D}_l}(R),\mathcal{D}'_l,\mathcal{D}'_r)$,
and in contrast with the case of the one-sided localization of $R$-modules,
Schelter proved that the subcategory $\mathcal{C}_{\sigma}$ of
complete triples is reflective in $\mathcal{W}_{\sigma}$. Our aim is
to extend the two-sided localization functor to bicategories (a
richer structure than a category) from which the study of Morita
equivalences of two-sided rings of quotients follows in a natural
way. Using these tools we obtain an exact sequence relating Picard groups of two-sided localizations,
generalizing substantially the well-known sequence obtained by Bass \cite{Bass} for the localization of
a commutative ring with respect to a multiplicative set of regular elements.

Now we summarize the contents of this paper.  In all this paper we
will basically follow notation and results from three sources,
Schelter \cite{Schelter2}, Goldman \cite{Goldman} and Stenstrom
\cite{Bo}. In section 1 given unital rings $S$ and $R$ and Gabriel
filters $\mathcal{H}_r$ and $\mathcal{D}_l$ of right ideals and left
ideals of $S$ and $R$, respectively, we will construct the two-sided
localization of an $R$-$S$-bimodule $M$ with respect to
$\mathcal{H}_r$ and $\mathcal{D}_l$, denoted by
$Q_{\mathcal{D}_l,\mathcal{H}_r}(M)$. In case that $M=R$ is the
regular $R$-$R$-bimodule then $Q_{\mathcal{D}_l, \mathcal{D}_r}(R)$
has a ring structure.

It is known that for unital rings $R$ and $S$ if the categories
Mod-$R$ and Mod-$S$ are equivalent then Mod-$\Qs(R)$ and
Mod-$\Qs(S)$ are equivalent categories too \cite{Ortega}. In the
general case our aim will be to determine what two-sided
localizations preserve Morita equivalence. It is well-known that
bicategories provide the best way to understand Morita equivalences
between rings and other structures, such as $C^*$-algebras, von
Neumann algebras, Lie groupoids \cite{Landsman}. Thus, in section 2
we are going to construct the bicategory $\mathcal{IB}_{\sigma}$
whose objects are the triples $(R,\mathcal{D}_l,\mathcal{D}_l)$
where $\mathcal{D}_r$ and $\mathcal{D}_l$ are Gabriel filters of
right and left ideals of $R$ respectively, and whose arrows between
$(R,\mathcal{D}_l,\mathcal{D}_l)$ and
$(S,\mathcal{H}_l,\mathcal{H}_l)$ are the equivalences between their
categories of modules that send the hereditary torsion theories
generated by $\mathcal{D}_r$ and $\mathcal{D}_l$ to the hereditary
torsion theories generated by $\mathcal{H}_r$ and $\mathcal{H}_l$
respectively. Then we will see that the two-sided localization
induces an endofunctor of the bicategory $\mathcal{IB}_{\sigma}$.
This fact will allow us to understand Morita equivalences in a
deeper way, identifying for example the bimodules that give rise to
the equivalences.

Finally, in section 3  we study the group of the auto-equivalences
of an object $(R,\mathcal{D}_{l},\mathcal{D}_r)$, written
$\textrm{Pic}(R,\mathcal{D}_{l},\mathcal{D}_r)$. This generalizes
the concept of $\textrm{Pic}(R)$, the Picard group of a ring $R$.
Indeed, we see that, in the case when $\mathcal{D}_r$ and
$\mathcal{D}_l$ are the Gabriel filters of dense right and left
ideals of $R$, respectively, we have
$\textrm{Pic}(R,\mathcal{D}_{l},\mathcal{D}_r)=\textrm{Pic}(R)$. For
a commutative ring $R$ and a multiplicative set of regular elements
$S$ Bass constructed an exact sequence that relates the Picard group
of $R$ and the Picard group of its localization $RS^{-1}$
\cite[Chapter III, Proposition 7.5]{Bass}. The machinery developed
in the previous chapters will allow us to extend this exact sequence
to a general two-sided localization of a non-commutative ring. Thus,
not only the commutative localization is extended to cover the case
of a two-sided Ore localization, but we also extend the type of
localization considered by Bass.

\section{Gabriel localizations and two-sided localizations.}

Throughout all rings will be associative and unital, and all modules
will be unitary. Given any ring $R$ let $(\mathcal{T},\mathcal{F})$
denote a hereditary torsion theory on the category Mod-$R$ of right
$R$-modules (see \cite[Chapter IV]{Bo}).

We recall (\cite[VI, Proposition 3.6 and Theorem 5.1]{Bo}) that
there is a bijection between hereditary torsion theories of Mod-$R$,
Gabriel filters of right ideals of $R$ and left exact radicals of
Mod-$R$. So we will talk of Gabriel filters $\Im$ rather than
hereditary torsion theories. Thus, given a Gabriel filter of right
(or left) ideals $\mathcal{D}$ of $R$ we denote by
$(\mathcal{T}_{\mathcal{D}},\mathcal{F}_{\mathcal{D}})$ the
hereditary torsion theory of $\textrm{Mod-}R$ (or $R\textrm{-Mod}$)
associated to $\mathcal{D}$.

\begin{exem}

\textbf{(1)} The hereditary torsion theory generated by the dense
right ideals $\Dr$ is of particular importance. $\Dr$ is the
strongest Gabriel filter for which $R$ is torsion-free.

\textbf{(2)} Let $S\subseteq R$ be a right Ore set of regular
elements, then the set of right ideals $\mathcal{D}^S=\{I\mid\
sR\subseteq I \text{ for some }s\in S\}$ is a Gabriel filter of
right ideals of $R$.

\end{exem}

Given a right $R$-module $M$ and a Gabriel filter of right ideals
$\mathcal{D}$, the \textit{module of quotients of  $M$ with respect
to $\mathcal{D}$}\index{module of quotients} is defined as
$$Q_{\mathcal{D}}(M):=\varinjlim_{I\in\mathcal{D}} \textrm{Hom}(I,\frac{M}{T_{\mathcal{D}}(M)})\,,$$
where $T_{\mathcal{D}}(M)=\{m\in M\mid mI=0 \text{ for some
}I\in\mathcal{D}\}$, the left exact radical associated to the
Gabriel filter $\mathcal{D}$. For every right $R$-module we have
that $Q_{\mathcal{D}}(M)$ has structure of right $R$-module, and if
$M=R_R$ then $Q_{\mathcal{D}}(R)$ has a ring structure \cite[Chapter
IX]{Bo}, and is called the \textit{right ring of quotients of $R$
with respect to $\mathcal{D}$}. Additionally to the right $R$-module
structure, $Q_{\mathcal{D}}(M)$ becomes a
$Q_{\mathcal{D}}(R)$-Module (see \cite[Chapter IX]{Bo}).

The most prominent examples of rings of quotients can be viewed as
localizations with respect to some Gabriel filters of right
ideals. Indeed, the maximal right ring of quotients is the
localization of $R$ with respect to $\Dr$ the filter of right
dense ideals of $R$, and for every right Ore set of regular
elements we have that $RS^{-1}$ is the localization of $R$ with
respect to the Gabriel filter of right ideals $\mathcal{D}^S$ of
$R$.

\begin{defi}
Let $R$ and $S$ be any rings, and $\mathcal{H}_r$ and
$\mathcal{D}_l$ Gabriel filters of right and left ideals of $S$
and $R$, respectively. Then we define
$_{\mathcal{D}_l}\Omega_{\mathcal{H}_r}$ as the set of right
ideals of $S\otimes R^{op}$ containing an ideal of the form
$H\otimes R^{op}+S\otimes D$ where $H\in \mathcal{H}_r $ and $D\in
\mathcal{D}_l$.
\end{defi}

If there is no confusion about which is the right and the left
Gabriel filter we denote $_{\mathcal{D}_l}\Omega_{\mathcal{H}_r}$
as $_l\Omega_r$. One can easily verify that $_l\Omega_r$ is a
Gabriel filter of right ideals of $S\otimes R^{op}$.

Thus given  an $R$-$S$-bimodule $_{R}M_S$, or equivalently, a right
$S\otimes R^{op}$-module, and $\mathcal{H}_r$ and $\mathcal{D}_l$
Gabriel filters of right and left ideals of $S$ and $R$,
respectively, we define the \textit{two-sided localization of $M$
with respect to $\mathcal{H}_r$ and $\mathcal{D}_l$} as
$$Q_{\mathcal{D}_l,\mathcal{H}_r }(M):=\varinjlim_{I\in
{_l\Omega_r}} \textrm{Hom}(I,\overline{M})\,,$$ where
$\overline{M}=\frac{M}{T_{_l\Omega_r}(M)}$.

A \textit{compatible pair $(f,g)$ on $M$} consists of an
$R$-homomorphism $f:D\rightarrow {}_{R}M$ and an $S$-homomorphism
$g:H\rightarrow M_S$ where $D\in \mathcal{D}_l$ and $H\in
\mathcal{H}_r$, which satisfy the compatibility condition
$xg(y)=(x)fy$ for every $x\in D$ and $y\in H$. We put $(f,g)\sim
(\overline{f},\overline{g})$ if and only if there exist $H'\in
\mathcal{H}_r$ and $D'\in \mathcal{D}_l$ such that
$f_{|D'}=\overline{f}_{|D'}$ and $g_{|H'}=\overline{g}_{|H'}$
Clearly $\sim$ is an equivalence relation on the set of compatible
pairs on $M$.

\begin{prop}\label{identificacio_compatibles}
Let $R$ and $S$ be rings and let $_{R}M_S$ be an $R$-$S$-bimodule.
Given Gabriel filters $\mathcal{D}_l$ and $\mathcal{H}_r$ of left
and right ideals of $R$ and $S$ respectively, there is a bijection
between elements of the two-sided localization
$Q_{\mathcal{D}_l,\mathcal{H}_r }(M)$ and the equivalence classes of
compatible pairs $(f,g)$ where $f: I\longrightarrow
{}_R\overline{M}$ and $g:J\longrightarrow \overline{M}_S$ for some
$I\in \mathcal{D}_l$ and $J\in \mathcal{H}_r$, with
$\overline{M}=\frac{M}{T_{_l\Omega_r}(M)}$.
\end{prop}
\begin{proof}
Let $h:H\otimes R^{op}+S\otimes D\longrightarrow \overline{M}$ be a
representative of an equivalence class of the two-sided localization
of $_{R}M_S$. Given any $x\in D$ and $y\in H$ we define
$(x)f:=h(1\otimes x)$ and $ g(y):=h(y\otimes 1)$. That $(f,g)$ is a
compatible pair follows in a straightforward way. Moreover, observe
that every representative of the equivalence class of $h$ gives rise
to a pair of compatible morphisms of the same equivalence class.

Conversely, given any compatible pair $(f,g)$ over
$_{R}\overline{M}_S$ we can construct the $S\otimes
R^{op}$-homomorphism $h:H\otimes R^{op}+S\otimes D\longrightarrow
\overline{M}$ such that $h(y\otimes 1)=g(y)$ and $h(1\otimes
x)=(x)f$ for every $y\in H$ and $x\in D$. This is well-defined by
the compatibility of the morphisms $f$ and $g$. Indeed, let us
construct the $S\otimes R^{op}$ morphisms
$$\begin{array}{cl} \alpha:  H\otimes R^{op} & \longrightarrow     \overline{M} \\
  y\otimes r & \longmapsto rg(y) \end{array}\qquad\text{and}\qquad
  \begin{array}{cl} \beta:  S\otimes D & \longrightarrow     \overline{M} \\
  s\otimes x & \longmapsto (x)fs \end{array}\,.$$

If $\sum y_i\otimes r_i=\sum s_j\otimes x_j$ we claim that
$\alpha(\sum y_i\otimes r_i)=\beta(\sum s_j\otimes x_j)$. Take
$s\in \cap s_j^{-1}H=I$. Then $s_js\in H$ for every $j$ and we
have that for every $r\in R$
$$\alpha(\sum y_i\otimes r_i)(s\otimes r)=\alpha((\sum y_i\otimes r_i)(s\otimes r))=\alpha((\sum s_j\otimes x_j)(s\otimes r))=$$
$$=\alpha(\sum s_js\otimes rx_j)=\sum rx_jg(s_js)=\sum (rx_j)fs_js=\beta(\sum s_j\otimes x_j)(s\otimes r)\,.$$

Then we have that $(\alpha(\sum y_i\otimes r_i)-\beta(\sum
s_j\otimes x_j))(I\otimes R^{op})=0$ where $I\in \mathcal{H}_r$.
Observe that symmetrically  we can construct $\cap Dr^{-1}_i=J\in
\mathcal{D}_l$ such that $(\alpha(\sum y_i\otimes r_i)-\beta(\sum
s_j\otimes x_j))(S\otimes J)=0$ and hence $(\alpha(\sum y_i\otimes
r_i)-\beta(\sum s_j\otimes x_j))(S\otimes J+ I\otimes R^{op})=0$.

But since $\overline{M}$ is $_l\Omega_r$-torsion-free, we get
$\alpha(\sum y_i\otimes r_i)-\beta(\sum s_j\otimes x_j)=0$ and
hence $\alpha(\sum y_i\otimes r_i)=\beta(\sum s_j\otimes x_j)$. So
$h$ is well-defined.
\end{proof}

From the above characterization of $Q_{\mathcal{D}_l,\mathcal{H}_r
}(M)$ the next lemma follows in a straightforward way
\begin{lem}
Let $R$ be a unital ring and let $\mathcal{D}_r$ and $\mathcal{D}_l$
be Gabriel filters of right and left ideals of $R$. Then
$Q_{\mathcal{D}_l,\mathcal{D}_r}(R)$ has a ring structure.
\end{lem}

Observe, that in the case that $M={}_{R}R_R$ with its Gabriel
filters of dense right and left ideals, $\Dr$ and $\Dl$
respectively, this coincides with the characterization of the
maximal symmetric ring of quotients of $R$ (see \cite{Lanning} and
\cite{Ortega}), i.e. ${Q}_{\Dl,\Dr}(R)=\Qs(R)$.

\section{Rings with filters of ideals as a bicategory.}

We can define a category  with objects the triples
$(R,\mathcal{D}_l,\mathcal{D}_r)$, where $R$ is a unital ring and
$\mathcal{D}_r$ and $\mathcal{D}_l$ are Gabriel filters of right and
left ideals of $R$, respectively, and the morphisms
$\alpha:(R,\mathcal{D}_l,\mathcal{D}_r)\longrightarrow
(S,\mathcal{H}_l,\mathcal{H}_r)$ are ring morphisms
$\alpha:R\longrightarrow S$ such that $\alpha(J)\cdot S\in
\mathcal{H}_r$ and $S\cdot\alpha(I)\in \mathcal{H}_l$ whenever $J\in
\mathcal{D}_r$ and $I\in \mathcal{D}_l$. It is proved in
\cite{Schelter} that the two-sided localization
$(R,\mathcal{D}_l,\mathcal{D}_r) \longrightarrow
(Q_{\mathcal{D}_l,\mathcal{D}_r}(R),\mathcal{D}'_l,\mathcal{D}'_r)$
defines an endofunctor of this category, where $\mathcal{D}'_l$ is
the Gabriel filter on $Q_{\mathcal{D}_l,\mathcal{D}_r}(R)$
consisting of left ideals of $Q_{\mathcal{D}_l,\mathcal{D}_r}(R)$
containing a left ideal of the form
$Q_{\mathcal{D}_l,\mathcal{D}_r}(R)\cdot I$, for $I\in
\mathcal{D}_l$, and similarly for $\mathcal{D}'_r$.

Now we define the oriented graph $\mathcal{B}_{\sigma}$ of unital
rings with Gabriel filters of right and left ideals whose objects
are the triples $\{(R,\mathcal{D}_l,\mathcal{D}_r)\mid R_R\in
\mathcal{F}_{\mathcal{D}_r}\text{ and }{}_RR\in
\mathcal{F}_{\mathcal{D}_l} \}$, and the morphisms between two
objects  $(R,\mathcal{D}_l,\mathcal{D}_r)$ and
$(S,\mathcal{H}_l,\mathcal{H}_r)$ are the $R$-$S$-bimodules
${}_{R}M_{S}$ such that:
\begin{enumerate}

\item[(Q1)] $M_S$ is a $\mathcal{H}_r$-torsion-free right $S$-module.

\item[(Q2)] ${}_RM$ is a $\mathcal{D}_l$-torsion-free left $R$-module.

\item[(Q3)] For every $D\in \mathcal{D}_r$ we have $(\frac{M}{DM})_S$ is
a $\mathcal{H}_r$-torsion module.

\item[(Q4)] For every $H\in \mathcal{H}_l$ we have ${}_R(\frac{M}{MH})$
is a $\mathcal{D}_l$-torsion module.

\end{enumerate}

The two-sided localization $Q:(R,\mathcal{D}_l,\mathcal{D}_r)
\longmapsto
(Q_{\mathcal{D}_l,\mathcal{D}_r}(R),\mathcal{D}'_l,\mathcal{D}'_r)$
defines mappings on the objects of $\mathcal{B}_{\sigma}$,
additionally we will see from the two following lemmas that $Q$
defines a morphism of graphs.

\begin{lem}\label{bimodule_localization}
Let  $(R,\mathcal{D}_l,\mathcal{D}_r)$ and
$(S,\mathcal{H}_l,\mathcal{H}_r)$ be two objects from
$\mathcal{B}_{\sigma}$ and let ${}_{R}M_S$ be a morphism between
them,  then $Q_{\mathcal{D}_l,\mathcal{H}_r}(M)$ is a
$Q_{\mathcal{D}_l,\mathcal{D}_r}(R)$-$Q_{\mathcal{H}_l,\mathcal{H}_r}(S)$-bimodule.
\end{lem}
\begin{proof}
Let $m\in Q_{\mathcal{D}_l,\mathcal{H}_r}(M)$, that is represented
by the pair of compatible homomorphisms $f:{}_{R}D\longrightarrow
{}_{R}M$, $g:{H}_{S}\longrightarrow {M}_{S}$, where $D\in
\mathcal{D}_l$ and $H\in \mathcal{H}_r$, and let $q\in
{Q}_{\mathcal{D}_l,\mathcal{D}_r}(R)$, that is represented by the
pair of compatible homomorphisms $f_1:{}_{R}I\longrightarrow R$,
$g_1:{J}_{R}\longrightarrow R$, where $I\in \mathcal{D}_l$ and $J\in
\mathcal{D}_r$.

Now taking the left ideal $\overline{D}:=I\cap (D)f_1^{-1}\in
\mathcal{D}_l $, we define the homomorphism
$\overline{f}:{}_R\overline{D}\rightarrow {}_RM$ by
$(x)\overline{f}:=((x)f_1)f$ for every $x\in \overline{D}$.

Now we define the following homomorphism
$$\begin{array}{cl} \widetilde{g_1}: JM_S   & \longrightarrow  M_S \\
  \sum y_i m_i  & \longmapsto \sum g_1(y_i) m_i  \end{array}
\,,$$ that is well-defined since ${}_{R}M$ is a
$\mathcal{D}_l$-torsion-free left $R$-module.

Note that by property $Q3$ we have that $(\frac{M}{JM})_S$ is a
$\mathcal{H}_r$-torsion module. So if we set the right ideal
$\overline{H}_S:=g^{-1}(JM)\in \mathcal{H}_r$ we define the
$S$-homomorphism $\overline{g}:\overline{H}_S\longrightarrow M_S$ as
$\overline{g}(y):=\widetilde{g_1}(g(y))$ for every $y\in
\overline{H}$.

Now we are going to see that the homomorphisms $\overline{f}$ and
$\overline{g}$ define a compatible pair
$(\overline{f},\overline{g})$. So given any $x\in \overline{D}$ and
$y\in \overline{H}$ we have $g(y)=\sum z_i m_i$ for some $z_i\in J$
and $m_i\in M$, so
$$x\overline{g}(y)=x\widetilde{g_1}(g(y))=x\widetilde{g_1}(\sum z_i
m_i)=x\sum g_1(z_i) m_i=$$
$$=(x)f_1\sum z_i
m_i=(x)f_1g(y)=((x)f_1)fy=(x)\overline{f}y\,.$$

Thus, the pair $(\overline{f},\overline{g})=:q\cdot m$ defines an
element of $Q_{\mathcal{D}_l,\mathcal{H}_r}(M)$, so we can define
the map $ {Q}_{\mathcal{D}_l,\mathcal{D}_r}(R)\times
Q_{\mathcal{D}_l,\mathcal{H}_r}(M) \longrightarrow
Q_{\mathcal{D}_l,\mathcal{H}_r}(M)$, $(q,m)\mapsto q\cdot m$, that
we can easily verify that induces  a structure of left
${Q}_{\mathcal{D}_l,\mathcal{D}_r}(R)$-module on
$Q_{\mathcal{D}_l,\mathcal{H}_r}(M)$.

Symmetrically we might see that $Q_{\mathcal{D}_l,\mathcal{H}_r}(M)$
is a right ${Q}_{\mathcal{H}_l,\mathcal{H}_r}(S)$-module, and that
is compatible with the left
${Q}_{\mathcal{D}_l,\mathcal{D}_r}(R)$-module structure.
\end{proof}

\begin{lem}\label{functor_weak_b}
Let  $(R,\mathcal{D}_l,\mathcal{D}_r)$ and
$(S,\mathcal{H}_l,\mathcal{H}_r)$ be two objects from
$\mathcal{B}_{\sigma}$ and let ${}_{R}M_S$ be a morphism between
them, then $Q_{\mathcal{D}_l,\mathcal{H}_r}(M)$ defines a morphism
between
 $({Q}_{\mathcal{D}_l,\mathcal{D}_r}(R),\mathcal{D}'_l,\mathcal{D}'_r)$ and $({Q}_{\mathcal{H}_l,\mathcal{H}_r}(S),\mathcal{H}'_l,\mathcal{H}'_r)$.
\end{lem}
\begin{proof}
We have seen above (Lemma \ref{bimodule_localization}) that
$Q_{\mathcal{D}_l,\mathcal{H}_r}(M)$ is a
${Q}_{\mathcal{D}_l,\mathcal{D}_r}(R)$-${Q}_{\mathcal{H}_l,\mathcal{H}_r}(S)$-bimodule,
so we have to check that $Q_{\mathcal{D}_l,\mathcal{H}_r}(M)$
verifies conditions $Q1-Q4$. First we are going to see that
$Q_{\mathcal{D}_l,\mathcal{H}_r}(M)$ is a
$\mathcal{H}'_r$-torsion-free right
${Q}_{\mathcal{H}_l,\mathcal{H}_r}(S)$-module. Indeed, let $m\in
Q_{\mathcal{D}_l,\mathcal{H}_r}(M)$ such that $mH'=0$ for some
$H'\in \mathcal{H}'_r$. Since $H'$ contains an ideal of the form
$H\cdot {Q}_{\mathcal{H}_l,\mathcal{H}_r}(S)$ for some $H\in
\mathcal{H}_r$ we have that $mH=0$. Now there exists a pair of
compatible homomorphisms $f:D\longrightarrow M$ and
$g:K\longrightarrow M$ where $D\in \mathcal{D}_l$ and $K\in
\mathcal{H}_r$ with $K\subseteq H$, such that $mk=g(k)$ and
$dm=(d)f$ for every $d\in D$ and $k\in K$. Thus, for every $d\in D$
and $k\in K$ we have that $0=dg(k)=(d)fk$, hence since $M$ is a
$\mathcal{H}_r$-torsion-free right $S$-module we get that $f=0$ and
$m=0$, as desired. Similarly we can prove $Q2$.

Now we shall see that $Q_{\mathcal{D}_l,\mathcal{H}_r}(M)$ verifies
$Q3$. Indeed, we must check that
$\frac{Q_{\mathcal{D}_l,\mathcal{H}_r}(M)}{D'\cdot
Q_{\mathcal{D}_l,\mathcal{H}_r}(M)}$ is a $\mathcal{H}'_r$-torsion
right ${Q}_{\mathcal{H}_l,\mathcal{H}_r}(S)$-module for every $D'\in
\mathcal{D}'_r$. Let $m\in Q_{\mathcal{D}_l,\mathcal{H}_r}(M)$,
since $\frac{Q_{\mathcal{D}_l,\mathcal{H}_r}(M)}{M}$ is a
$_{l}\Omega_r$-torsion $S\otimes R^{op}$-module there exist $I\in
\mathcal{D}_l$ and $J\in \mathcal{H}_r$ with $Im,mJ\subseteq M$. Now
let $D\in \mathcal{D}_r$ with $D\cdot
{Q}_{\mathcal{D}_l,\mathcal{D}_r}(R)\subseteq D'$, since $M$
satisfies $Q3$ for every $y\in J$ there exists $H_y\in
\mathcal{H}_r$ such that $m y H_y\subseteq DM$. Now we set
$K:=\sum_{y\in J} yH_y$ that belongs to $\mathcal{H}_r$ by a
property of the Gabriel filters (namely property T4 in \cite[page
146]{Bo}). Thus, we have that $m\cdot(K \cdot
{Q}_{\mathcal{D}_l,\mathcal{D}_r}(R)) \subseteq D\cdot
Q_{\mathcal{D}_l,\mathcal{H}_r}(M)\subseteq D'\cdot
Q_{\mathcal{D}_l,\mathcal{H}_r}(M)$, as desired. Symmetrically we
can check $Q4$.
\end{proof}

Thus, we have that
$$\begin{array}{cl} Q:\mathcal{B}_{\sigma} &  \longrightarrow\mathcal{B}_{\sigma} \\
(R,\mathcal{D}_l,\mathcal{D}_r)  & \longrightarrow  ({Q}_{\mathcal{D}_l,\mathcal{D}_r}(R),\mathcal{D}'_l,\mathcal{D}'_r)\\
  M  & \longmapsto  Q_{\mathcal{D}_l,\mathcal{H}_r}(M)   \end{array}
\,,$$ defines a functor on $\mathcal{B}_{\sigma}$.

We would like to equip the graph $\mathcal{B}_{\sigma}$ of rings
with Gabriel filters of right and left ideals with a richer
structure than a category. So we recall the well-known definition of
bicategory (\cite{Benabou} and \cite{Borceaux}).

\begin{defi}
A \textit{bicategory} $\mathcal{B}$ consists of the following
structures:
\begin{enumerate}
\item A class of objects $A$, $B$, $C$, ... called 0-cells;

\item For each pair of 0-cells $A$ and $B$, a small category
$\mathcal{B}(A,B)$ which objects $\mathcal{B}^1(A,B)$ are called
1-cells, and which morphisms are called 2-cells.

\item For each triple $A$, $B$ and $C$ of objects (0-cells), there is
a composition law given by a functor
$\otimes_{A,B,C}:\mathcal{B}(A,B)\times\mathcal{B}(B,C)\longrightarrow
\mathcal{B}(A,C)$.

\item For each object $A$, there is an identity object functor
$\textrm{I}_A:\textbf{1}\longrightarrow\mathcal{B}(A,A)$, where
$\textbf{1}$ stands for the final object in the category of small
categories, such that:
\begin{enumerate}
\item For every 0-cells $A,B,C,D$  we have natural associativity
isomorphisms
$$((-\otimes_{A,B,C}-)\otimes_{A,C,D}-) \simeq
(-\otimes_{A,B,D}(-\otimes_{B,C,D}-))\,.$$

\item For every $A$ and $B$ we have natural unity isomorphisms
$(-\otimes \textrm{I}_B)\simeq (\textrm{I}_A\otimes -)\simeq
\textrm{Id}_{\mathcal{B}(A,B)}$.

\end{enumerate}
The above isomorphisms must satisfy some coherence properties (for
details see \cite[7.7]{Borceaux}).
\end{enumerate}
\end{defi}

\begin{exem} Our starting point will be the bicategory of rings and
bimodules $\mathcal{B}im$:

(0-cells) Rings $R$, $S$,...

(1-cells) Given $R$ and $S$ then $\mathcal{B}(R,S)= R$-$S$-Bimod.

(2-cells) Given ${}_{R}M_S$ and ${}_{R}N_S$ the 2-cells are the
$R$-$S$-bimodule homomorphisms.

The composition law is given by the tensor product, and given any
ring $R$ its identity object is the regular $R$-$R$-bimodule
${}_{R}R_R$. \end{exem}

One could be tempted to define a bicategorical structure on
$\mathcal{B}_{\sigma}$ using tensor product as composition of
$1$-cells. Unfortunately, conditions $Q1$-$Q4$ are not stable under
tensor product. To correct this default, we have to restrict our
attention to a smaller class of bimodules.

It is known  \cite[Chapter II]{Bass} that every equivalence $T:
\textrm{Mod-}R\longrightarrow \textrm{Mod-}S$ is determined by an
invertible $R$-$S$-bimodule $_{R}M_{S}$, i.e. $T\simeq -\otimes_R
M=:T_M$ (see \cite{Anderson_fuller} for further information about
invertible bimodules). Thus, there exists a Morita context $\alpha:
M\otimes M^* \longrightarrow R$ and $\beta: M^*\otimes_R M
\longrightarrow S$, with $\alpha$ and $\beta$ associative bimodule
isomorphisms (see \cite[Chapter II, Proposition 3.1]{Bass}), where
$M^*=\textrm{Hom}_R(M,R)$ is the dual $R$-module which is indeed an
$S$-$R$-bimodule. Observe that the equivalence of the categories of
left modules is given by $U_{M}:=M\otimes_S -$. Let us denote by
$\mathcal{IB}im(R,S)$ the class of the invertible $R$-$S$-bimodules,
that corresponds to the invertible 1-cells of the bicategory
$\mathcal{B}im$.

It is known \cite[Chapter X \S 3]{Bo} that every natural equivalence
$T_M$ induces a bijective correspondence between the Gabriel filters
of right ideals of $R$ and the Gabriel filters of right ideals of
$S$, denoted by $T_M^{\sharp}$. Observe that for every invertible
$R$-$S$-bimodule $M$ and Gabriel filters $\mathcal{D}_r$ and
$\mathcal{H}_r$ of right ideals of $R$ and $S$, respectively, we
have that $T^{\sharp}_M(\mathcal{D}_r)=\mathcal{H}_r$ is equivalent
to $T_M(\mathcal{T}_{\mathcal{D}_r})=\mathcal{T}_{\mathcal{H}_r}$
and $T_M(\mathcal{F}_{\mathcal{D}_r})=\mathcal{F}_{\mathcal{H}_r}$.

Now we define the bicategory $\mathcal{IB}_{\sigma}$ whose 0-cells
are the triples $\{(R,\mathcal{D}_l,\mathcal{D}_r)\mid R_R\in
\mathcal{F}_{\mathcal{D}_r}\text{ and }{}_RR\in
\mathcal{F}_{\mathcal{D}_l} \}$, and given two 0-cells
$(R,\mathcal{D}_l,\mathcal{D}_r)$ and
$(S,\mathcal{H}_l,\mathcal{H}_r)$ we set the category
$\mathcal{IB}_{\sigma}((R,\mathcal{D}_l,\mathcal{D}_r),(S,\mathcal{H}_l,\mathcal{H}_r))$
which objects (1-cells) are the $R$-$S$-bimodules $M$ that satisfy
\begin{enumerate}

\item[(Q)] ${}_{R}M_S$ is an invertible $R$-$S$-bimodule such that
$T^{\sharp}_M(\mathcal{D}_r)=\mathcal{H}_r$ and
$U^{\sharp}_{M}(\mathcal{H}_l)=\mathcal{D}_l$, \end{enumerate} and
with the $R$-$S$-bimodules homomorphism as morphisms (2-cells).
Observe that bimodules that satisfy condition $(Q)$ automatically
satisfy conditions $(Q1)$-$(Q4)$.

In this case, it is easy to check that the tensor product of
bimodules acts as composition law in the 1-cells of
$\mathcal{IB}_{\sigma}$, thus $\mathcal{IB}_{\sigma}$ is a
bicategory.

Now we would like to prove that the mapping $Q$ restricted to
$\mathcal{IB}_{\sigma}$ induces a functor of bicategories.

\begin{rema}\label{unit}
For every invertible $R$-$S$-bimodule $M$ there exist some
$\widetilde{m}_1,\ldots,\widetilde{m}_r\in M$ and
$\widetilde{n}_1,\ldots,\widetilde{n}_r\in M^*$ such that
$1_R=\sum_j \alpha(\widetilde{m}_j\otimes\widetilde{n}_j)$.
\end{rema}

\begin{prop}\label{iso_loc_tensor} Let
$(R,\mathcal{D}_l,\mathcal{D}_r)$ and
$(S,\mathcal{H}_l,\mathcal{H}_r)$ be two $0$-cells of
$\mathcal{IB}_{\sigma}$, then for every 1-cell ${}_RM_S$ between
them, the maps $$\begin{array}{cl} \varphi: M\otimes_S {Q}_{\mathcal{H}_l,\mathcal{H}_r}(S)   & \longrightarrow  Q_{\mathcal{D}_l,\mathcal{H}_r}(M) \\
  \sum m_i\otimes q_i  & \longmapsto  \sum m_i  q_i
  \end{array}\qquad\text{and}\qquad
\begin{array}{cl} \phi: {Q}_{\mathcal{D}_l,\mathcal{D}_r}(R)\otimes_RM & \longrightarrow  Q_{\mathcal{D}_l,\mathcal{H}_r}(M) \\
  \sum p_j\otimes m_j  & \longmapsto  \sum p_j m_j  \end{array}  \,.$$
 are
$R$-$S$-bimodule isomorphisms.
\end{prop}
\begin{proof}
Let $\alpha: M\otimes M^* \longrightarrow R$ and $\beta: M^*\otimes
M \longrightarrow S$ be a Morita context, with $\alpha$ and $\beta$
bimodule isomorphisms. First we are going to prove that $\varphi$ is
injective. Suppose that there exists $\sum m_i\otimes q_i\in
M\otimes_S {Q}_{\mathcal{H}_l,\mathcal{H}_r}(S)$ such that
$\varphi(\sum m_i\otimes q_i)=\sum m_i q_i=0$, hence using Remark
\ref{unit} we have
$$\sum_i
m_i\otimes q_i=\sum_i 1_R\cdot m_i\otimes q_i=\sum_i (\sum_j
\alpha(\widetilde{m}_j\otimes\widetilde{n}_j)m_i)\otimes q_i=$$
$$=\sum_i(\sum_j \widetilde{m}_j \beta(\widetilde{n}_j\otimes m_i))\otimes q_i=\sum_j(\widetilde{m}_j\otimes \sum_i\beta(\widetilde{n}_j\otimes m_i)q_i)\,.$$

There exists $H\in \mathcal{H}_r$ such that $q_iH\subseteq S$ for
every $i$, and $\sum_i\beta(\widetilde{n}_j\otimes
m_i)q_ih=\beta(\widetilde{n}_j\otimes \sum_i
 m_iq_i h)=0$ for every $h\in H$. But since
 ${Q}_{\mathcal{H}_l,\mathcal{H}_r}(S)$ is
 $\mathcal{H}_r$-torsion-free it yields that
 $\sum_i\beta(\widetilde{n}_j\otimes m_i)q_i=0$ and hence $\sum_i
m_i\otimes q_i=0$. Thus, $\varphi$ is injective.

Now we are going to see that $\varphi$ is also surjective. Indeed,
first we will prove that we can extend $\beta$ to $\overline{\beta}:
M^* \otimes_R Q_{\mathcal{D}_l,\mathcal{H}_r}(M) \longrightarrow
{Q}_{\mathcal{H}_l,\mathcal{H}_r}(S)$ as follows. Let us consider
any $n\in N$ and $m\in Q_{\mathcal{D}_l,\mathcal{H}_r}(M)$, that is
represented by the compatible pair $f:D\rightarrow {}_{R}M$ and
$g:H\rightarrow M_S$ where $D\in \mathcal{D}_l$ and $H\in
\mathcal{H}_r$. We define the $S$-homomorphism $\overline{g}:H
\longrightarrow S$ as $\overline{g}(h):=\beta(n\otimes g(h))$. Now,
we define the following $S$-homomorphism
$$\begin{array}{cl} \widetilde{f}: {}_SM^*\cdot D   & \longrightarrow S \\
  \sum n_i d_i  & \longmapsto  \sum \beta(n_i \otimes (d_i)f)  \end{array}
\,.$$

We have that $\widetilde{f}$ is well-defined since $S$ is
$\mathcal{H}_r$-torsion-free. Now if we consider the
$S$-homomorphism $\textrm{R}_n:S\longrightarrow M^*$, that is right
multiplication by $n$, we have that $K:=\textrm{R}_n^{-1}(M^*D)\in
\mathcal{H}_l$ since ${}_{S}(\frac{M^*}{M^*\cdot D})$ is a
$\mathcal{H}_l$-torsion module. So we can define the
$S$-homomorphism $\overline{f}:K \longrightarrow S$ as
$(k)\overline{f}:=(kn)\widetilde{f}$. Now we shall see that the pair
$(\overline{f},\overline{g})$ is a compatible one. Indeed, let $k\in
K$ and $h\in H$, then
$$(k)\overline{f}h=(kn)\widetilde{f}h=(\sum n_i d_i)\widetilde{f}h=\sum\beta(n_i\otimes(d_i)f)h=\sum\beta(n_i\otimes(d_i)fh)=$$
$$=\sum\beta(n_i\otimes d_ig(h))=\beta(\sum n_i d_i \otimes g(h))=\beta(kn \otimes g(h))=k\overline{g}(h)\,.$$

Thus, the pair $(\overline{f},\overline{g})$ defines an element of
${Q}_{\mathcal{H}_l,\mathcal{H}_r}(S)$. Moreover it is clear that
$\overline{\beta}_{|M^*\otimes M}=\beta$, and that we have the
following compatibility $\alpha(m \otimes
n)\overline{m}=m\overline{\beta}(n \otimes \overline{m})$ for every
$n\in N$, $m\in M$ and $\overline{m}\in
Q_{\mathcal{D}_l,\mathcal{H}_r}(M)$.

Now using that we can extend $\beta$ to $\overline{\beta}$ we shall
prove that $\varphi$ is surjective. Indeed, let $\overline{m}$ be
any element from $Q_{\mathcal{D}_l,\mathcal{H}_r}(M)$, we have that
$\overline{m}=1_R \overline{m}=\sum_i \alpha(\widetilde{m}_i
\otimes\widetilde{n}_i) \overline{m}=\sum \widetilde{m}_i
\overline{\beta}(\widetilde{n}_i\otimes\overline{m})$. So it follows
that $\varphi(\sum_i \widetilde{m}_i\otimes
\overline{\beta}(\widetilde{n}_i
\otimes\overline{m}))=\overline{m}$, so we have that $\varphi$ is a
surjective homomorphism. Thus, we get $M\otimes_S
{Q}_{\mathcal{H}_l,\mathcal{H}_r}(S)\cong
Q_{\mathcal{D}_l,\mathcal{H}_r}(M)$. A symmetric argument gives the
other isomorphism.
\end{proof}

Now with this setting, given the $0$-cells
$(R,\mathcal{D}_l,\mathcal{D}_r)$, $(S,\mathcal{H}_l,\mathcal{H}_r)$
and $(T,\mathcal{K}_l,\mathcal{K}_r)$ from $\mathcal{IB}_{\sigma}$
and any 1-cells ${}_{R}M_S\in
\mathcal{IB}_{\sigma}^{1}((R,\mathcal{D}_l,\mathcal{D}_r),(S,\mathcal{H}_l,\mathcal{H}_r))$
and ${}_{S}N_T\in
\mathcal{IB}_{\sigma}^{1}((S,\mathcal{H}_l,\mathcal{H}_r),(T,\mathcal{K}_l,\mathcal{K}_r))$
it follows from Proposition \ref{iso_loc_tensor} that we have the
following commutative diagram of isomorphisms:

$$\xymatrix{ { (M\otimes_S N)\otimes_T Q(T)}\ar[r]^{\varphi}_{\cong}\ar[d]^{\cong} & {{Q}_{\mathcal{D}_l,\mathcal{K}_r}(M\otimes_S N)}\ar[dd]_{h}^{\cong}  &  {{Q}(R)\otimes_R (M\otimes_S N)}\ar[l]_{\phi}^{\cong}\ar[d]^{\cong} \\
{M\otimes_S {Q}_{\mathcal{H}_l,\mathcal{K}_r}(N)}\ar[d]^{\cong} & &
{{Q}_{\mathcal{D}_l,\mathcal{H}_r}(M) \otimes_S N}\ar[d]^{\cong} \\
{M\otimes_S {Q}(S)\otimes_{Q(S)}
{Q}_{\mathcal{H}_l,\mathcal{K}_r}(N)}\ar[r]^{\cong} &
{{Q}_{\mathcal{D}_l,\mathcal{H}_r}(M)\otimes_{Q(S)}
{Q}_{\mathcal{H}_l,\mathcal{K}_r}(N)} &
{{Q}_{\mathcal{D}_l,\mathcal{H}_r}(M) \otimes_{Q(S)} {Q}(S)\otimes_S
N}\ar[l]_{\cong}}$$

\begin{corol}\label{Q_functor}
$Q$ defines a functor in the bicategory $\mathcal{IB}_{\sigma}$.
Therefore, if $T_M:\textrm{Mod-}R\longrightarrow \textrm{Mod-}S$ is
an equivalence then ${Q}_{\mathcal{D}_l,\mathcal{D}_r}(R)$ and
${Q}_{U_{M^*}^{\sharp}(\mathcal{D}_l),T_M^{\sharp}(\mathcal{D}_r)}(S)$
are Morita equivalent rings.
\end{corol}
\begin{proof}
We define the functor $Q$ by
$$\begin{array}{cl} Q: \mathcal{IB}_{\sigma}& \longrightarrow \mathcal{IB}_{\sigma} \\
  (R,\mathcal{D}_l,\mathcal{D}_r)  & \longmapsto  ({Q}_{\mathcal{D}_l,\mathcal{D}_r}(R),\mathcal{D}'_l,\mathcal{D}'_r) \\
  {}_RM_S & \longmapsto  {Q}_{\mathcal{D}_l,\mathcal{H}_r}(M)  \\
(f:M\rightarrow N)  & \longmapsto
(\overline{f}:Q_{\mathcal{D}_l,\mathcal{H}_r}(M)\rightarrow
Q_{\mathcal{D}_l,\mathcal{H}_r}(N))  \end{array} \,,$$ for $M\in
\mathcal{IB}^1_{\sigma}((R,\mathcal{D}_l,\mathcal{D}_r),(S,\mathcal{H}_l,\mathcal{H}_r))$,
where $\mathcal{D}'_r=\{I_{{Q}(R)}\mid D\cdot
{Q}_{\mathcal{D}_l,\mathcal{D}_r}(R)\subseteq I\text{ for some }
D\in \mathcal{D}_r \}$ and $\mathcal{D}'_l=\{{}_{{Q}(R)}J\mid
{Q}_{\mathcal{D}_l,\mathcal{D}_r}(R) \cdot D \subseteq J\text{ for
some } D\in \mathcal{D}_l \}$. The mapping between the 1-cells is
well-defined by Lemma \ref{functor_weak_b} and Proposition
\ref{iso_loc_tensor}, and the above commutative diagram shows that
$Q$ commutes with the tensor product.
\end{proof}

We have that the Gabriel filter of the dense right ideals of $R$
is the biggest filter with respect to which the regular right
$R$-module $R_R$ is torsion-free, thus since an equivalence of
modules induces a bijection between the torsion theories of both
categories the following result seems natural:

\begin{lem}\label{equivalence_dense}(c.f.\cite[Chapter X, Proposition 3.2]{Bo})
Let $T:\textrm{Mod-}R\longrightarrow \textrm{Mod-}S$ be an
equivalence, and let $\mathcal{I}_{dr}(R)$ be the Gabriel filter
of dense right ideals of $R$. Then
$T^{\sharp}(\mathcal{I}_{dr}(R))=\mathcal{I}_{dr}(S)$ is the
Gabriel filter of dense right ideals of $S$.
\end{lem}
\begin{corol}\label{symmetric_morita_equivalence}
If  $R$ and $S$ are two Morita equivalent rings, then $\Qs(R)$ and
$\Qs(S)$ are Morita equivalent rings too.
\end{corol}

\section{The Picard group.}

Given a bicategory $\mathcal{B}$ we define its \textit{large
Picard groupoid}\index{Picard!groupoid}, written as
$\textrm{Pic}(\mathcal{B})$, as the groupoid of isomorphism
classes of  invertible arrows of $\mathcal{B}$ (see \cite{Benabou}
for details). In the bicategory $\mathcal{B}$im of rings and
bimodules the invertible arrows are the invertible bimodules. So
$\textrm{Pic}(\mathcal{B}\textrm{im})$ is the groupoid of all the
isomorphism classes $[T]$ of natural equivalences
$T:R\text{-Mod}\longrightarrow S\text{-Mod}$ and the
multiplication is given by composition of the functors, or
equivalently by the tensor product of invertible bimodules.

The orbits of $\textrm{Pic}(\mathcal{B}\textrm{im})$ are the Morita
equivalence classes, and given a ring $R$ the isotropy group of the
autoequivalences of $R$-Mod is the \textit{Picard
group}\index{Picard!group} of $R$, written as $\textrm{Pic}(R)$. The
group law is induced by the composition of functors. As we saw
before, we can see $\textrm{Pic}(R)$ as the group of isomorphisms
classes $[{}_{R}M_R]$ of invertible $R$-$R$-bimodules. The group law
is given by $[M][N]=[M\otimes N]$, and with $[M]^{-1}=[M^*]$.

Let us consider the bicategory $\mathcal{IB}_{\sigma}$ defined
previously, then given a $0$-cell $(R,\mathcal{D}_l,\mathcal{D}_r)$
we denote its isotropy group as
$\textrm{Pic}(R,\mathcal{D}_l,\mathcal{D}_r)$. Observe that by Lemma
\ref{equivalence_dense} we have that
$\textrm{Pic}(R,\mathcal{D}_l,\mathcal{D}_r)=\textrm{Pic}(R)$ when
$\mathcal{D}_r$ and $\mathcal{D}_l$ are the filters of the dense
right and left ideals respectively.

Moreover, by Corollary \ref{Q_functor} we have the following group
homomorphism
$$\begin{array}{cl} Q: \textrm{Pic}(R,\mathcal{D}_{l},\mathcal{D}_r) & \longrightarrow  \textrm{Pic}({Q}_{\mathcal{D}_l,\mathcal{D}_r}(R),\mathcal{D}'_{l},\mathcal{D}'_r)\\
 {[M]} & \longmapsto {[{Q}_{\mathcal{D}_l,\mathcal{D}_r}(M)]}\end{array}\,.$$

Observe that in particular we have the following group
homomorphism
$$\begin{array}{cl} Q: \textrm{Pic}(R) & \longrightarrow  \textrm{Pic}(\Qs(R))\\
 {[M]} & \longmapsto  {[M\otimes_R \Qs(R)]=[\Qs(R)\otimes_R M]} \end{array}\,.$$

Now given $(R,\mathcal{D}_l,\mathcal{D}_r)$ and any
$R$-$R$-subbimodule $M$ of ${Q}_{\mathcal{D}_l,\mathcal{D}_r}(R)$ we
define $M^{-1}:=\{n\in {Q}_{\mathcal{D}_l,\mathcal{D}_r}(R)\mid
nM,Mn\subseteq R\}$. Observe that $M^{-1}$ is also a
$R$-$R$-bimodule. Moreover we can define natural associative
$R$-$R$-bimodule homomorphisms $\alpha: M\otimes_R M^{-1}
\longrightarrow  R $ and $ \beta: M^{-1}\otimes_R M \longrightarrow
R$.

Thus, we define the group $\textrm{Pic}(R\mid
{Q}_{\mathcal{D}_l,\mathcal{D}_r}(R))$ of $R$-$R$-subbimodules $M$
of ${Q}_{\mathcal{D}_l,\mathcal{D}_r}(R)$ such that
$M^{-1}M=M^{-1}M=R$ ($M$ \textit{invertible in}
${Q}_{\mathcal{D}_l,\mathcal{D}_r}(R)$ ), and such that
$T^{\sharp}_M(\mathcal{D}_r)=\mathcal{D}_r$ and
$U^{\sharp}_{M}(\mathcal{D}_l)=\mathcal{D}_l$. Observe that
$\textrm{Pic}(R\mid {Q}_{\mathcal{D}_l,\mathcal{D}_r}(R))$ is a
group with unit $R$, with operation given by the product $M_1\cdot
M_2$ that by the following result is isomorphic to $M_1\otimes_R
M_2$.
\begin{prop} Let $M$ and $N$ be $R$-$R$-bimodules of
${Q}_{\mathcal{D}_l,\mathcal{D}_r}(R)$ with $M\in\textrm{Pic}(R\mid
{Q}_{\mathcal{D}_l,\mathcal{D}_r}(R))$. Then $M\otimes_R N\cong MN$.
\end{prop}
\begin{proof}
It is a straightforward exercice that
$$\begin{array}{cl} f: M\otimes_R N & \longrightarrow  MN \\
 \sum m_i\otimes n_i & \longmapsto  {\sum m_i n_i} \end{array}\,$$
is a $R$-$R$-bimodule isomorphism.\end{proof}

\begin{rema}
If $\Dl$ and $\Dr$ are the Gabriel filters of dense right and left
ideals of $R$ respectively, then by Lemma \ref{equivalence_dense},
$\textrm{Pic}(R\mid {Q}_{\Dl,\Dr}(R))=\textrm{Pic}(R\mid \Qs(R))$ is
the set of  $R$-$R$-subbimodules $M$ of $\Qs(R)$ with
$MM^{-1}=M^{-1}M=R$.
\end{rema}

Given a ring $R$ we denote by $\mathcal{U}(\mathcal{Z}(R))$ the set
of all central invertible elements of $R$.

\begin{theor}\label{exact_local}
Let $(R,\mathcal{D}_l,\mathcal{D}_r)$ be a $0$-cell in
$\mathcal{IB}_{\sigma}$. Then the sequence
$$0\longrightarrow \mathcal{U}(\mathcal{Z}(R))\longrightarrow\mathcal{U}(\mathcal{Z}({Q}_{\mathcal{D}_l,\mathcal{D}_r}(R)))\longrightarrow \textrm{Pic}(R\mid{Q}_{\mathcal{D}_l,\mathcal{D}_r}(R))\longrightarrow $$
$$\longrightarrow\textrm{Pic}(R,\mathcal{D}_l,\mathcal{D}_r) \longrightarrow \textrm{Pic}({Q}_{\mathcal{D}_l,\mathcal{D}_r}(R),\mathcal{D}'_l,\mathcal{D}'_r)$$
is exact.
\end{theor}
\begin{proof}
It is easy to check that $\mathcal{Z}(R)\subseteq
\mathcal{Z}({Q}_{\mathcal{D}_l,\mathcal{D}_r}(R))$ using that $R$ is
torsion-free with respect to $\mathcal{D}_r$ and $\mathcal{D}_l$. It
is clear that the kernel of
$$\begin{array}{cl} \varphi_1: \mathcal{U}(\mathcal{Z}({Q}_{\mathcal{D}_l,\mathcal{D}_r}(R))) & \longrightarrow  \textrm{Pic}(R\mid {Q}_{\mathcal{D}_l,\mathcal{D}_r}(R))\\
 x & \longmapsto  xR \end{array}$$
is $\mathcal{U}(\mathcal{Z}(R))$.

Now, we claim that the kernel of
$$\begin{array}{cl} \varphi_2: \textrm{Pic}(R\mid {Q}_{\mathcal{D}_l,\mathcal{D}_r}(R)) & \longrightarrow  \textrm{Pic}(R,\mathcal{D}_l,\mathcal{D}_r)\\
 M & \longmapsto  {[M]} \end{array}$$
consists of the $R$-$R$-subbimodules of
${Q}_{\mathcal{D}_l,\mathcal{D}_r}(R)$ of the form $xR$ with $x\in
\mathcal{U}(\mathcal{Z}({Q}_{\mathcal{D}_l,\mathcal{D}_r}(R)))$.
Indeed, $[M]=[R]$ implies that there exists a $R$-$R$-bimodule
isomorphism $f:R\longrightarrow M$, and in this case $f$ can be
represented by multiplication by an element $x:=f(1)$. For every
$a\in R$ we have that $xa=f(1)a=f(a)=af(1)=ax$ so $x$ commutes with
every element of $R$ and hence $x\in
\mathcal{Z}({Q}_{\mathcal{D}_l,\mathcal{D}_r}(R))$ and we get
$M=xR$. Since $M$ is invertible and $x\in
\mathcal{Z}({Q}_{\mathcal{D}_l,\mathcal{D}_r}(R))$ there is $y\in
M^{-1}$ such that $xy=1$, so $x\in
\mathcal{U}(\mathcal{Z}({Q}_{\mathcal{D}_l,\mathcal{D}_r}(R)))$.

Now, we are going to see that the kernel of
$$\begin{array}{cl} \varphi_3: \textrm{Pic}(R,\mathcal{D}_l,\mathcal{D}_r) & \longrightarrow  \textrm{Pic}({Q}_{\mathcal{D}_l,\mathcal{D}_r}(R),\mathcal{D}'_l,\mathcal{D}'_r)\\
 {[M]}  & \longmapsto  {[M\otimes {Q}_{\mathcal{D}_l,\mathcal{D}_r}(R)]=[{Q}_{\mathcal{D}_l,\mathcal{D}_r}(M)]} \end{array}$$
are those classes $[M]$ of $R$-$R$-bimodules $M$ isomorphic  to some
$\overline{M}\in \textrm{Pic}(R\mid
{Q}_{\mathcal{D}_l,\mathcal{D}_r}(R))$.

Let $M$ be a $R$-$R$-bimodule such that $\varphi:
{Q}_{\mathcal{D}_l,\mathcal{D}_r}(M) \longrightarrow
{Q}_{\mathcal{D}_l,\mathcal{D}_r}(R)$ is a
${Q}_{\mathcal{D}_l,\mathcal{D}_r}(R)$-${Q}_{\mathcal{D}_l,\mathcal{D}_r}(R)$-bimodule
isomorphism. Since $M$ is a $_l\Omega_r$-torsion-free
$R$-$R$-bimodule we have that $M$ is a $R$-$R$-subbimodule of
${Q}_{\mathcal{D}_l,\mathcal{D}_r}(M)$. Thus, $M$ is isomorphic to
the $R$-$R$-subbimodule $\overline{M}:=\varphi(M)$ of
${Q}_{\mathcal{D}_l,\mathcal{D}_r}(R)$. We have to see that
$\overline{M}$ is invertible in
${Q}_{\mathcal{D}_l,\mathcal{D}_r}(R)$. Indeed, let $N$ be the
inverse $R$-$R$-bimodule of $M$, and let $\alpha:M\otimes_R N
\longrightarrow R$ and $\beta:N\otimes_R M \longrightarrow R$ be
$R$-$R$-bimodule isomorphisms that satisfy the associativity
property. Using isomorphisms of Proposition \ref{iso_loc_tensor} we
can extend $\alpha$ and $\beta$ to isomorphisms
$$\xymatrix{ {\overline{\alpha}:{Q}_{\mathcal{D}_l,\mathcal{D}_r}(M)\otimes_R
N}\ar[r]^{\cong} &  {{Q}_{\mathcal{D}_l,\mathcal{D}_r}(R)\otimes_R
M\otimes_R N}\ar[rr]^{\textrm{Id}\otimes\alpha} & &
{{Q}_{\mathcal{D}_l,\mathcal{D}_r}(R)\otimes_R R\cong
{Q}_{\mathcal{D}_l,\mathcal{D}_r}(R)}}\,,$$ and
$$\xymatrix{ {\overline{\beta}: N\otimes_R{Q}_{\mathcal{D}_l,\mathcal{D}_r}(M)}\ar[r]^{\cong}  &
{N\otimes_R
M\otimes_R{Q}_{\mathcal{D}_l,\mathcal{D}_r}(R)}\ar[rr]^{\beta\otimes\textrm{Id}}
& & {R\otimes_R {Q}_{\mathcal{D}_l,\mathcal{D}_r}(R)\cong
{Q}_{\mathcal{D}_l,\mathcal{D}_r}(R)}}\,.$$

Then applying the functor $N\otimes_R -$  and then
$\overline{\beta}$ to the homomorphism
$$\xymatrix{ {R}\ar@{^{(}->}[r]^{i} & {{Q}_{\mathcal{D}_l,\mathcal{D}_r}(R)}\ar[r]^{\varphi^{-1}}_{\cong} & {{Q}_{\mathcal{D}_l,\mathcal{D}_r}(M)}}\,, $$
we get the homomorphism $\phi_1$ given by the composition
$$\xymatrix{ {N}\ar[r]^{\cong} & {N\otimes_RR}\ar@{^{(}->}[rr]^{\textrm{Id}\otimes i} & & {N\otimes_R{Q}_{\mathcal{D}_l,\mathcal{D}_r}(R)}\ar[rr]^{\textrm{Id}\otimes \varphi^{-1}}_{\cong} & & {N\otimes_R{Q}_{\mathcal{D}_l,\mathcal{D}_r}(M)}\ar[r]^{\overline{\beta}} &  {{Q}_{\mathcal{D}_l,\mathcal{D}_r}(R)} }\,. $$

Observe that since $N_R$ is a finitely generated projective right
$R$-module we have that $\phi_1$ is a $R$-$R$-bimodule monomorphism,
thus we have that $N$ is isomorphic to the $R$-$R$-subbimodule
$N_1:=\phi_1(N)$ of ${Q}_{\mathcal{D}_l,\mathcal{D}_r}(R)$.
Symmetrically applying the functor $-\otimes_R N$ and
$\overline{\alpha}$ we construct the $R$-$R$-bimodule monomorphism
$\phi_2$, so we get that $N$ is isomorphic to the
$R$-$R$-subbimodule $N_2:=\phi_2(N)$ of
${Q}_{\mathcal{D}_l,\mathcal{D}_r}(R)$. Now we claim that
$N_1\overline{M}=R$. Indeed, let $\overline{m}\in \overline{M}$ and
$\overline{n}\in N_1$, we have that there exist $m\in M$ and $n\in
N$ such that $\varphi(m)=\overline{m}$ and $\phi_1(n)=\overline{n}$.
Now we have that
$$\overline{n}\cdot\overline{m}=\phi_1(n)\varphi(m)=\overline{\beta}(n\otimes\varphi^{-1}(1))\varphi(m)=\overline{\beta}(n\otimes\varphi^{-1}(1)\varphi(m))=$$
$$=\overline{\beta}(n\otimes\varphi^{-1}(\varphi(m)))=\overline{\beta}(n\otimes m)=\beta(n\otimes m)\in R\,.$$

Now, since there exist $n_i\in N$ and $m_i\in M$ such that $\sum
\beta(n_i\otimes m_i)=1$ we have that $\sum
\phi_1(n_i)\varphi(m_i)=1$ and hence $N_1\overline{M}=R$ as desired.
A symmetric argument shows that $\overline{M}N_2=R$.  Thus, we have
that $N_1=N_1R=N_1(\overline{M}N_2)=(N_1\overline{M})N_2=RN_2=N_2$
by the associativity of the ring
${Q}_{\mathcal{D}_l,\mathcal{D}_r}(R)$. So $\overline{M}$ is
invertible in ${Q}_{\mathcal{D}_l,\mathcal{D}_r}(R)$ as desired.

On the other hand it is clear that if $M\in \textrm{Pic}(R\mid
{Q}_{\mathcal{D}_l,\mathcal{D}_r}(R))$ then $M\otimes_R
{Q}_{\mathcal{D}_l,\mathcal{D}_r}(R)\cong
{Q}_{\mathcal{D}_l,\mathcal{D}_r}(R)$. Indeed, we can define the
$R$-${Q}_{\mathcal{D}_l,\mathcal{D}_r}(R)$-bimodule isomorphism
$$\begin{array}{cl} \varphi:M\otimes_R
{Q}_{\mathcal{D}_l,\mathcal{D}_r}(R) & \longrightarrow{Q}_{\mathcal{D}_l,\mathcal{D}_r}(R)\\
\sum m_i\otimes q_i  & \longmapsto  \sum m_i q_i
\end{array}\,.$$

\end{proof}

\begin{corol}
Let $R$ be a unital ring, then the sequence
$$0\longrightarrow \mathcal{U}(\mathcal{Z}(R))\longrightarrow\mathcal{U}(\mathcal{Z}(\Qs(R)))\longrightarrow \textrm{Pic}(R\mid\Qs(R))\longrightarrow \textrm{Pic}(R) \longrightarrow \textrm{Pic}(\Qs(R))$$
is exact.
\end{corol}

Now we are going to study what happens with the Picard groups of a
two-sided Ore localization. Let $S$ be a two-sided Ore set of
regular elements of $R$. Let us consider the ring of fractions
$RS^{-1}=S^{-1}R$, that is the localization of $R$ with respect to
$\mathcal{D}^{S}_r$ and $\mathcal{D}^{S}_l$ the Gabriel filters of
right and left ideals of $R$ that have nonempty intersection with
$S$ respectively.

\begin{corol}
Let $R$ be a unital ring and let $S$ be a two-sided Ore set of
regular elements of $R$, then the sequence
$$0\longrightarrow \mathcal{U}(\mathcal{Z}(R))\longrightarrow\mathcal{U}(\mathcal{Z}(RS^{-1}))\longrightarrow \textrm{Pic}(R\mid RS^{-1})\longrightarrow \textrm{Pic}(R,\mathcal{D}^{S}_r,\mathcal{D}^{S}_l) \longrightarrow \textrm{Pic}(RS^{-1})$$
is exact, where $\textrm{Pic}(R\mid RS^{-1})$ is the set of all the
$R$-$R$-subbimodules $M$ of  $RS^{-1}$ such that
$MM^{-1}=M^{-1}M=R$.
\end{corol}
\begin{proof}
We apply Theorem \ref{exact_local} to get the desired exact
sequence. We claim that for every invertible $R$-$R$-subbimodule $M$
in $RS^{-1}$ we have that
$T^{\sharp}_M(\mathcal{D}_r^{S})=\mathcal{D}_r^{S}$ and
$U^{\sharp}_M(\mathcal{D}_l^{S})=\mathcal{D}_l^{S}$. Indeed, we only
have to check that $T_M(\frac{R}{sR})=\frac{R}{sR}\otimes_R M\cong
\frac{M}{sM}$ is a $\mathcal{D}_r^{S}$-torsion right $R$-module for
every $s\in S$. So let $s\in S$, for every $m\in M$ we have to find
$t\in S$ with $mt\in sM$. We can write $m=au^{-1}$ with $a\in R$ and
$u\in S$, and there exist $b\in R$ and $w\in S$ such that $aw=sb$,
so
$$m(uw)=(au^{-1})(uw)=aw=sb\,,$$ with $uw\in S$, but $b\in R$.
Since $MM^{-1}=R$ there exist $m_i\in M$ and $n_i\in M^{-1}$ for
$i\in \{1,\ldots,k\}$ with $b=\sum m_in_i$. Now we can write
$n_i=c_iz_i^{-1}$ with $c_i\in R$ and $z_i\in S$ for every $i\in
\{1,\ldots,k\}$. We can find $z\in S$ and $\overline{c}_{i}\in R$
such that $n_i=\overline{c}_{i} z^{-1}$ for every $i\in
\{1,\ldots,k\}$. So $bz=\sum m_i n_i z= \sum m_i
(\overline{c}_{i}z^{-1}z)=\sum m_i \overline{c}_{i} \in M$, and
hence
$$m(uwz)=sbz=s(\sum m_i \overline{c}_{i})$$ with $uwz\in S$ and $\sum m_i \overline{c}_{i}\in M$, as desired.

Thus, $T_M(\frac{R}{sR})$ is a $\mathcal{D}_r^{S}$-torsion right
$R$-module for every $s\in S$, so
$T^{\sharp}_M(\mathcal{D}_r^{S})\subseteq\mathcal{D}_r^{S}$, and
applying the same argument with $M^{-1}$ in place of $M$ we get the
reverse inclusion, so that
$T^{\sharp}_M(\mathcal{D}_r^{S})=\mathcal{D}_r^{S}$. Symmetrically
we have that $U^{\sharp}_M(\mathcal{D}_l^{S})=\mathcal{D}_l^{S}$.

Finally, observe that since the induced filters are
$(\mathcal{D}^{S}_r)'=\{RS^{-1}\}$ and
$(\mathcal{D}^{S}_r)'=\{RS^{-1}\}$ we have that
$\textrm{Pic}(RS^{-1},(\mathcal{D}^{S}_r)',(\mathcal{D}^{S}_l)')=\textrm{Pic}(RS^{-1})$.
\end{proof}

We have seen above that $\textrm{Pic}(R\mid RS^{-1})$ coincides with
all the invertible $R$-$R$-subbimodules of $RS^{-1}$. However, one
can construct rings $R$ and choose non-maximal two-sided Ore sets
$S\subseteq R$ such that
$\textrm{Pic}(R,\mathcal{D}^{S}_r,\mathcal{D}^{S}_l)$ does not
coincide with all the isomorphism classes of invertible
$R$-$R$-bimodules $\textrm{Pic}(R)$. \newline

\textbf{Open question:} When $S$ is the maximal two-sided Ore set of
regular elements of $R$, determine whether
$\textrm{Pic}(R,\mathcal{D}^{S}_r,\mathcal{D}^{S}_l)$ coincides with
$\textrm{Pic}(R)$.
\newline

Observe that an affirmative answer to this question would imply that
the maximal two-sided Ore localization of two unital Morita
equivalent rings are Morita equivalent rings.

Finally, we will compare the sequence of Theorem \ref{exact_local}
with the one obtained by Bass in \cite[Chapter III, Proposition
7.5]{Bass} for commutative rings. Let $R$ be a unital commutative
ring. Given any right $R$-module $M$ and ring automorphisms $\phi$
and $\varphi$ we denote by ${}_{\varphi}M_{\phi}$ the
$R$-$R$-bimodule whose additive group is $M$ and whose bimodule
structure is given by $r\cdot m=m\varphi(r)$ and $ m\cdot
r=m\phi(r)$ for every $m\in M$ and $r\in R$.

We denote by $\textrm{Pic}_R(R)$ the group of the isomorphism
classes of the invertible right $R$-modules (see \cite{Bass}).

It is known \cite[Chapter II \S 5]{Bass} that given any
$R$-$R$-bimodule $M$ the action of $R$ on the left is determined
by an automorphism $\varphi$ of $R$, that is $r\cdot
m=m\varphi(r)$ for every $r\in R$ and $m\in M$. So it follows that
given any $R$-$R$-bimodule $M$ there exists a right $R$-module $N$
and a ring automorphism $\varphi$ such that $M\cong
{}_{\varphi}N_1$ as $R$-$R$-bimodules. Thus we get the following
result.
\begin{prop}\cite[Chapter II, Proposition 5.4]{Bass}
Let $R$ be a commutative ring, then
$$0\longrightarrow \textrm{Pic}_R(R)\longrightarrow \textrm{Pic}(R)\longrightarrow Aut(R)\longrightarrow 0$$
is exact, and splits by $\varphi\mapsto [{}_{\varphi}R_{1}]$.
\end{prop}

Thus, we can see every element of $\textrm{Pic}(R)$ as a pair
$([M],\varphi)=[{}_{\varphi}M_1]=[{}_{\varphi}R_1\otimes_R
{}_1M_1]$ where $M$ is an invertible right $R$-module and
$\varphi$ is a ring automorphism. Then given $([M],\varphi)$ and
$([N],\phi)$ we have that ${}_{\varphi}M_1\otimes_R
{}_{\phi}N_1\cong {}_{\phi\circ\varphi}R_{1}\otimes_R
{}_{1}(M\otimes_R {}_{\phi}N_1)_1\in ([M\otimes_R
{}_{\phi}N_1],[\phi\circ\varphi])$.

Thus, we have that $\textrm{Pic}(R)\cong \textrm{Pic}_R(R)\rtimes
\textrm{Aut}(R)$ with the above product.

Now, let $S$ be a saturated multiplicative set of regular elements
of $R$, and let $RS^{-1}$ be its ring of fractions. Recall that
$RS^{-1}$ is the localization of $R$ with respect to the Gabriel
filter of ideals of $R$ that have nonempty intersection with $S$,
denoted by $\mathcal{D}^{S}$.

Now, we are going to describe $\textrm{Pic}(R,\mathcal{D}^{S})$,
that is, the group of isomorphism classes $[M]$ of invertible
$R$-$R$-bimodules such that
$T^{\sharp}_M(\mathcal{D}^{S})=\mathcal{D}^{S}$. Since every
invertible $R$-$R$-bimodule is isomorphic to ${}_{\varphi}M_1$ for
some right $R$-module $M$ and ring automorphism $\varphi$, we get
that $T^{\sharp}_M(\mathcal{D}^{S})=\mathcal{D}^S$ if and only if
$\varphi(S)=S$.

The group of the ring automorphisms such that $\varphi(S)=S$ is
denoted by $\textrm{Aut}^{S}(R)$.
\begin{lem}
Let $R$ be a unital commutative ring, and let $S$ be a saturated
multiplicative system of regular elements of $R$, then
$$\textrm{Pic}(R,\mathcal{D}^{S})\cong \textrm{Pic}_R(R)\rtimes
\textrm{Aut}^S(R)\,.$$
\end{lem}

Now, observe that if $S$ is a saturated multiplicative system of
regular elements of $R$, and we consider the ring of fractions of
$R$ with respect to $S$, by Theorem \ref{exact_local} we get the
following exact sequence
$$\vcenter{\xymatrix{{0}\ar[r] & {\mathcal{U}(R)}\ar[r] & {\mathcal{U}(RS^{-1})}\ar[r] & {\textrm{Pic}(R\mid RS^{-1})}\ar[r] &  {} }}$$
$$\vcenter{\xymatrix{{}\ar[r] &  {\textrm{Pic}_R(R)\rtimes
\textrm{Aut}^S(R)}\ar[r] & {\textrm{Pic}_{RS^{-1}}(RS^{-1})\rtimes
\textrm{Aut}(RS^{-1})}}}\,,$$ where

$$\begin{array}{cl} \varphi_2: \textrm{Pic}(R\mid RS^{-1}) & \longrightarrow  \textrm{Pic}_R(R)\rtimes
\textrm{Aut}^S(R) \\
 {M}  & \longmapsto  {([M],\textrm{Id})} \end{array}\,,$$
and
$$\begin{array}{cl} \varphi_3: \textrm{Pic}_R(R)\rtimes
\textrm{Aut}^S(R) & \longrightarrow
\textrm{Pic}_{RS^{-1}}(RS^{-1})\rtimes
\textrm{Aut}(RS^{-1}) \\
 {([M],\varphi)}  & \longmapsto  {([MS^{-1}],\overline{\varphi})} \end{array}\,,$$
where $\overline{\varphi}: RS^{-1}\longrightarrow RS^{-1}$ is the
ring automorphism defined as
$\overline{\varphi}(ab^{-1})=\varphi(a)\varphi(b)^{-1}$, that is
well-defined since $\varphi\in \textrm{Aut}^S(R)$. Splitting off the
Aut-terms, we obtain:

\begin{corol}\cite[Chapter III, Proposition 7.5]{Bass}
Let $R$ be a unital commutative ring, and let $S$ be a saturated
multiplicative set of nonzero divisors of $R$. Then the sequence
$$0\longrightarrow \mathcal{U}(R)\longrightarrow\mathcal{U}(RS^{-1})\longrightarrow \textrm{Pic}(R\mid RS^{-1})\longrightarrow \textrm{Pic}_{R}(R) \longrightarrow \textrm{Pic}_{RS^{-1}}(RS^{-1})$$
is exact.
\end{corol}

\section*{Acknowledgments}

This research was partially supported by MEC-DGESIC (Spain) through
Project MTM2005-00934, and by the Comissionat per Universitats i
Recerca de la Generalitat de Catalunya. Part of this paper was
written while the author visited Queen's University Belfast as part
of a research project funded by Royal Society. I gratefully
acknowledge my advisor P.Ara for all his useful comments and help
that led to the presentation of this paper.

\end{document}